# Influence tests I: ideal composite hypothesis tests, and causal semimeasures

Bruno Bauwens *

November 5, 2018

### Abstract

Ratios of universal enumerable semimeasures corresponding to hypotheses are investigated as a solution for statistical composite hypotheses testing if an unbounded amount of computation time can be assumed.

Influence testing for discrete timeseries is defined using generalized structural equations. Several ideal tests are introduced, and it is argued that when Halting information is transmitted, in some cases, instantaneous cause and consequence can be inferred where this is not possible classically.

The approach is contrasted with Bayesian definitions of influence, where it is left open whether all Bayesian causal associations of universal semimeasures are equal within a constant. Finally the approach is also contrasted with existing engineering procedures for influence and theoretical definitions of causation.

## Introduction

The paper introduces necessary tools to define and investigate general purpose influence tests for discrete timeseries in the ideal case that an unbounded amount of computation time is available.

In statistics, a *simple hypothesis* corresponds to a hypothesis that contains enough information to infer a unique probability distribution over all measurement data that could be expected. Comparing two simple hypothesis according to given data is well understood [22]. In many cases too little useful information is available to construct such a unique probability distribution. There is no general accepted solution for the problem of composite hypotheses testing. It is argued here that the ratio test for universal semimeasures can theoretically define a solution for composite hypothesis testing.

An extensive literature exists on the definition of 'influence' both in statistics and in philosophy [7, 17, 21, 25, 31]. However, most of this work defines influence only when a fixed distribution is already available. General purpose statistical tests, differ from simple statistical tests with respect that they can infer a useful notion of influence without reference to a semimeasure. To define ideal influence tests, or to interpret practical algorithms such as [10, 15, 16, 24, 26, 28, 30, 32], there is no theory available that considers both the statistical interpretation and the computability aspects.

---

*Department of Electrical Energy, Systems and Automation, Ghent University, Technologiepark 913, B-9052, Ghent, Belgium, Bruno.Bauwens@ugent.be. Supported by a Ph.D grant of the Institute for the Promotion of Innovation through Science and Technology in Flanders (IWT-Vlaanderen).

Causality is often related to structural equations [25]. Traditionally, computable functions are used to study these generalized structural equations. The set of semimeasures corresponding to these structural equations do not lead to sets of semimeasures with a universal element. To solve this problem, structural equations with partial computable functions will be considered.

**Overview of results.** In section 1, composite hypothesis testing using universal semimeasures is discussed. It is shown that if a composite hypothesis is associated with a set of semimeasures that is "testable" and is a product of convex sets of semimeasures, it has a universal enumerable semimeasure among its enumerable semimeasures. This result will be applied to the hypothesis of independence and of timeseries being influence-free.

In section 2, different hypothesis of influence-free and causal timeseries are defined. All corresponding sets of semimeasures have universal elements and there values can differ substantially if Halting information is instantaneously transmitted between two timeseries. This happens depending on whether instantaneous information is assumed to be originated from a hidden source, or is instantaneously transmitted from the first signal to the second or opposite.

In section 3, causal semimeasures Bayasianally associated with enumerable semimeasure are introduced, which are a superset of the enumerable causal semimeasures, defined in section 2. It is shown that they do not have a universal element. However, for the causal semimeasures associated with some universal element, it is left open whether such a universal element exists. The results are summarized in figure 3.

The logarithm of the proposed ideal statistical tests, define an algorithmic variant of the Shannon information transfer. In section 4 both quantities are related, and therefore an alternate interpretation of the algorithms in [28, 32, 26, 24, 32] can be given. Also Granger causality can be interpreted in this framework. Finally, the proposed test for influence is contrasted to Shannon information transfer of the minimal sufficient statistic, and to graphical representations of minimal sufficient statistics as in [18]. It is shown that enumerable algorithmic information transfer determines plausible causal relations, where these relations can not be determined from probabilistic minimal sufficient statistics.

**Definitions and notation.** Let $\omega$ be the set of natural numbers, and $\epsilon$ be the empty string. The binary strings $2^{<\omega}$ of finite length can be associated with $\omega$. Let $l(x)$ denote the length of $x$ in its binary expansion and let $2^n$ be the set of strings $x$ with $l(x) = n$. Let $\omega^n$ be the set of $\omega$-sequences of length $n$. Let $\Phi$ be a universal Turing machine. $\Phi_t(p, x) \downarrow = y$ means that $\Phi$ on input $p, x$ outputs $y$, and halts in less than $t$ computation steps. Let $2^\omega$ denote Cantor space, the space of infinite binary sequences. $2^\omega$ can be associated with the Real numbers in $[0, 1]$.[1] For $r \in [0, 1]$, $r^k$ are the first $k$ decimals in the binary expansion of $r$. For $x \in 2^{<\omega}$, $x^k$ denotes $x_1 x_2 ... x_k$. A real function $f : \omega \to [0, 1]$ is computable if there is a string $p$ such that for all $k, x$: $\Phi(p, x, k) \downarrow = f(x)^k$. An enumeration of an one-argument real function $f(x)$ is a two-argument computable rational function $g(x, t)$ such that for all $t$: $g(u, t) \leqslant g(u, t+1)$ and such that $\lim_t g(u, t) = f(u)$. With abuse of notation an enumeration of $f$ is denoted as $f_t$. A function $f$ is enumerable, respectively co-enumerable, iff $f$, respectively $-f$ has an enumeration. A semimeasure $P$ is a non-negative real function that satisfies $\sum \{P(x) : x \in \omega\} \leqslant 1$. A semimeasure $P$ (multiplicatively) dominates a semimeasure $Q$, notation: $P \geqslant^* Q$, if a constant $c$ exists such that for all $x$: $cP(x) \geqslant Q(x)$ [20]. $P =^* Q$, if $P \leqslant^* Q$ and $Q \leqslant^* P$. A set $S$ of semimeasures has a universal element $m$ if $m \in S$ and $m$ dominates all semimeasures in $S$. A function $f : \omega \to \omega$ (additively) dominates a function $g : \omega \to \omega$, notation: $f \geqslant^+ g$, iff there is a constant $c$ such that for all $x$: $f(x) + c \geqslant g(x)$. $f =^+ g$ iff $f \leqslant^+ g$ and $g \leqslant^+ f$. For some of the proofs prefix-free Kolmogorov complexity is needed, which is defined in the appendix.

---

[1] In this association for example the number 0.5 is both associated with 0111... as with 1000....

# 1 Ideal Bayesian composite hypothesis testing

This section describes scientific hypothesis testing and defines a notion of significance for the ratio test of universal semimeasures. Proposition 1.2 shows the existence of universal semimeasures for a variety of hypotheses.

## 1.1 What do hypothesis tests?

Scientific modeling is the process of making *rules* for symbols that represent *observables* or properties of observables. These rules can be iteratively applied and combined to reproduce past observations or predict and control future observations in specific *contexts*. In scientific modeling, one often first starts to infer rules from a restricted context (inference), and conjectures that they apply to more contexts (generalization).

Scientists agree or disagree on the applicability of rules and models in different contexts. When a rule is under discussion, the rule is called a *hypothesis*. If scientists agree on the applicability of hypotheses and models, science is advancing. When probabilities are involved in a model, a hypothesis can imply a semimeasure over all possible expected observations in a context of an experiment. A hypothesis that implies such a semimeasure is called a *simple hypothesis*. The discussion of the applicability of the hypothesis in the restricted context of the experiment often happens through the use of significance or hypothesis tests applied to observations of an experiment which is called data. Such tests define significance or probability of type I error, and sensitivity, or probability of type II error, or one minus power. The test rejects or fails to reject the zero hypothesis if the significance is below or above a predefined value. The test favors the alternate hypothesis if also sensitivity is below a predefined value.

In scientific models, probabilities arise when in a context some variables are not observable or some variables are beyond control. Three types of probabilities can be obtained, either

- from repetitive observations of data. (Frequentistic probabilities)

- from rules about the variables that are not subject to discussion. (Objectivistic probabilities)

- from an unknown observer-dependent model. (Subjectivistic probabilities)

Significance and sensitivity are probabilities. If the probabilities of the hypothesis are obtained frequentistically or objectively, then at least significance and sensitivity can have the same interpretation. For a frequentistic setting, the significance respectively sensitivity of a statistical test is the maximal limit fraction of repetitive evaluations of the test where $P^0$ is rejected, in a context where $P^0$ describes the observed data well, respectively will not reject $P^0$ in a context where $P^A$ describes the observed data well. In the objectivistic setting, the significance respectively sensitivity is the objective prior probability that the zero hypothesis disqualifies itself, respectively alternate hypothesis disqualifies itself.

More formal, let an one-sided statistical test $d(x)$ be high for data $x$ that seem to contradict the zero hypothesis. Assume that zero and alternate hypothesis are simple with probability distributions $P^0$ and $P^A$. Significance and sensitivity for data $x$ according to $d(x)$ are given by:

$$\alpha(x) = \sum \{P^0(y) : d(y) \geqslant d(x)\}$$

$$\beta(x) = \sum \{P^A(y) : d(y) \leqslant d(x)\}.$$

If $\alpha(x)$ is small, a scientist will conclude that either a rare event has occurred or that $P^0$ is not representative for $x$. In practice, he will reject the zero hypothesis. If also $\beta(x)$ is small, he will favor the alternate hypothesis.

A specific choice of such a statistical test $d$ is given by the likelihood ratio $P^A(x)/P^0(x)$. This ratio has also an interpretation within Bayesian statistics: if $a^0/a^1$ represents the ratio of prior belief in the hypothesis corresponding to $P^0$ relative to the belief in the hypothesis corresponding to $P^A$, than after observing data $x$ the posterior ratio of the beliefs is:

$$\frac{a^0}{a^1} \frac{P^0(x)}{P^A(x)}.$$

Due to the coding theorem [20] this Bayesian interpretation is also justified by an Occam's razor argument that favors the hypothesis that can be described with minimal code length.

The Newman-Pearson lemma states that $\beta \circ \alpha^{-1} : [0,1] \rightarrow [0,1]$ is uniformly maximal for

$$d(x) = P^A(x)/P^0(x).$$

This means that there is a test that has for any significance an optimal sensitivity. This shows that optimal hypothesis testing is equivalent to likelihood ratio testing. Remark that significance and sensitivity are bounded by $P^0(x)/P^A(x)$.

A *composite hypothesis* is a collection of rules that imply a set of semimeasures. There is no accepted optimal general way of extending the hypothesis testing from simple hypothesis testing to composite hypothesis testing. Many methods in literature are proposed that are theoretical optimal under some conditions [9], or have been found to be empirically useful in specific contexts. Let $H^0$ and $H^A$ be the sets of semimeasures constituted by the zero and alternate hypothesis.

- *Uniformly optimal test:* In specific cases, there is a test that has an optimal $\beta \circ \alpha^{-1}$ function for all combination of tests in $H^0$ and $H^A$.

- *Bayesian approaches*: Assign some fixed prior probability to all semimeasures, this reduces the problem to simple hypothesis testing. Often it is not possible to extend the hypothesis with an acceptable prior and therefore this is a subjective method.

- *Generalized maximal likelihood*: This is the likelihood ratio of the best case hypotheses:
  $$\frac{\max\{P(x) : P \in H^A\}}{\max\{P(x) : P \in H^0\}}.$$
  This is the most commonly used method. In specific cases this method is proved to be optimal, but in other cases it has problems or is subject to discussion [9].

Suppose that a composite zero $H^0$ and alternate hypothesis $H^A$ have universal semimeasures $m^0$ and $m^1$. If the sets of the semimeasures are convex and countable $H^i = P_1^i, P_2^i, \ldots$ for $i = 0, A$, the universal semimeasures $m^i$ satisfy:

$$m^i =^* \sum_j a_j P_j^i. \tag{1}$$

This means that hypothesis selection by the likelihood ratio

$$d(x) = m^A(x)/m^0(x),$$

can be considered as a Bayesian approach to composite hypothesis selection. Since all universal semimeasures are equal up to a constant factor, the subjectivity is limited to a constant factor. Also because $m^i \in H^i$ and because it is multiplicatively optimal, it can be considered as generalized maximal likelihood testing if one can neglect the constant factors.

The significance of $d$ relative to $m^0$, does not have a direct frequentistic or objective interpretation. In general a repetitive experiment with controlled and uncontrolled variables in the environment can not frequentistically evaluate to $m^0$ or $m^A$. Objectively, neither $m^0$ or $m^A$ are guaranteed to become the accepted semimeasures for the context under discussion. Assume $H^0$ and $H^A$ satisfy the conditions of Proposition 1.2. Since any choice of the positive constants $a_i$ in (1) results in a universal semimeasure, without loss of generality the constants $a_i = 1/i(\log i)^2$ can be chosen. Let $i \leqslant k/(\log k)^3$ than:

$$\frac{P_i^0(x)}{m^A(x)} \leqslant^* \frac{km^0(x)}{m^A(x)}.$$

If the significance of $d$ is large, this can mean that:

- Some complex model from the zero hypothesis describes the data.

- The alternate hypothesis $m^A$ better describes the data.

- A rare event has occurred.

In many cases the first interpretation must be partly taken into account, and therefore one should look for a separate notion of significance for the statistic $d(x)$. For example, a frequentistic significance bound is obtained by a permutation test for the Shannon information transfer statistic in [24].

## 1.2 Universal semimeasures for a composite hypothesis.

Semimeasures are used in stead of measures, since the set of computable or enumerable measures has no universal element. The set of co-enumerable semimeasures also has no universal element [4, 20].

A positive real function $P$ is a *length conditional semimeasure*, if for all $n$:

$$\sum \{P(x) : x \in 2^n\} \leqslant 1.$$

The use of length conditional semimeasures allows to reduce technical details. Furthermore they can be justified by remarking that in many experimental setups, the amount of generated data, is fixed before the experiment starts. From now on, *semimeasure* is short for *length conditional semimeasure*.

**Definition 1.1.** Let $S$ be a set of semimeasures:

- $S^\uparrow$ is the subset of enumerable semimeasures in $S$.

- $S$ is *testable* iff there is a computable logic expression $L$ such that for any semimeasure $P$: $P \in S$ iff some rational approximation $P_t$ of $P$ satisfies:

$$\forall t, n \leqslant t : L(P_t^n),$$

where $P_t^n$ is the finite restriction of $P_t$ on $2^n$.

- $S$ is *convex* iff from any $P, Q \in S$, and $a, b \in [0, 1]$ with $a + b \leqslant 1$: $aP + bQ \in R$.

- The product set of two sets of semimeasures $S, T$ is given by

$$S \times T = \{PQ : P \in S \wedge Q \in T\}$$

Remark that the product set of two semimeasures is also a semimeasure.

**Proposition 1.2.** *Let $S, T$ be sets of semimeasures.*

(i) *If $S$ is testable and contains $P_0 = 0$ than $S^\uparrow$ is enumerable.*

(ii) *If $S$ is convex and $S^\uparrow$ can be enumerated as $P_1, P_2, ...,$ than $S^\uparrow$ contains the universal semimeasure*

$$m^{S^\uparrow} = \sum a_i P_i,$$

*where $a_i > 0$ is any computable real function such that $\sum_{i \in \omega} a_i \leqslant 1$.*

(iii) *If $S^\uparrow, T^\uparrow$ have universal elements $m^S, m^T$, then*

$$m^{S^\uparrow \times T^\uparrow} = m^S m^T$$

*is a universal element for $S^\uparrow \times T^\uparrow$.*

*Proof.* The first two items of the proposition are a direct generalization of the proof of the existence of universal enumerable semimeasures [12, 19, 20].
**Part (i).** Define the enumeration $P_{i,t}$: let $P_{i,0}(x) = 0$ for all $i, x$. Remark that $P_{i,0} \in S$. For all $t$ let

$$P_{i,t}(x) = \max\{\Phi_t(i, x, s) : s \leqslant t \wedge \Phi_t(i, x, s) \downarrow\},$$

if $\sum\{P_{i,t}(x) : x \in 2^n\} \leqslant 1$, and $L(P_{i,t}^n)$ is true, otherwise let for all $x$:

$$P_{i,t}(x) = P_{i,t-1}(x).$$

Remark that for all $i, t$, $P_{i,t}(x)$ is computable and that if $i$ is a code for a $Q \in S^\uparrow$, than there is an $i$ such that $Q = P_i$.
**Part (ii).** Let:

$$m_t^S = \sum\{a_i P_{i,t} : i \leqslant t\}.$$

Remark that $m_t^S$ is computable, that it increases with $t$, and therefore $m^S$ is enumerable. Remark that by convexity for all $t$, $m_t^S \in S^\uparrow$, and for any $n$, the values $m_t^S(w)$ with $w \in 2^{<n}$ remain constant for some $t$ large enough. Therefore the limit is also in $S^\uparrow$. Finally remark that $m^S$ dominates all $P_i$.
**Part (iii).** Clearly $m^{S^\uparrow \times T^\uparrow} \in S^\uparrow \times T^\uparrow$. Let $R \in S^\uparrow \times T^\uparrow$. It remains to show that $R \leqslant^* m^{S^\uparrow \times T^\uparrow}$. There exist $P \in S^\uparrow, Q \in T^\uparrow$ such that $R = PQ$. since $P \leqslant c_P m^S$ and $Q \leqslant c_Q m^T$, we have that

$$R = PQ \leqslant c_P c_Q m^S m^T = c_P c_Q m^{S^\uparrow \times T^\uparrow}.$$

$\square$

From Proposition 1.2 it follows that the set of univariate, bivariate and conditional enumerable semimeasures have a universal element denoted as: $m(x)$, $m(x, y)$, and $m(x|y)$. The set of independent enumerable semimeasures is given by $P(x, y) = Q(x)R(y)$, for $Q, R$ univariate semimeasures. The set satisfies the conditions of item ($iii$) of Proposition 1.2, and therefore has universal element $m(x)m(y)$. Also remark that by Corollary 2.6 there are sets $S, T$, such that the universal element of $S^\uparrow \times T^\uparrow$ can be a factor $o(n/\log n)$ lower than the universal element of $(S \times T)^\uparrow$.

By the coding theorem, we have that $-\log m(x) =^+ K(x)$ and $-\log m(x, y) =^+ K(x, y)$. For this reason $I(x; y) = K(x) + K(y) - K(x, y)$ naturally appears as some notion of confidence for ideal independence [19].

## 2 Influence-free and causal semimeasures

This section derives several influence tests from generalized structural equations both for pairs of discrete variables, and for pairs of discrete timeseries. It is shown that when Halting information is present in two observations $x, y$, the obtained universal elements from the structural equation hypotheses can imply slightly different likelihoods if $x$ is assumed to cause $y$ or $y$ is assumed to cause $x$. When $x, y$ represent discrete timeseries, the difference in likelihood can become significant depending on whether $x$ is assumed to instantaneously cause $y$ or $y$ is assumed to instantaneously cause $x$.

### 2.1 Statistical explanatory model

First the concept of statistical explanatory model, is discussed within the computability framework.

Cantor space $2^\omega = [0, 1]$ with tree topology is assumed, it is, for every $r \in 2^{<\omega}$:

$$[r] = \{\alpha \in 2^{<\omega} : r \sqsubset \alpha\},$$

with $r \sqsubset \alpha$ meaning that $r$ is a prefix of $\alpha$. The measure is given by $\mu([r]) = 2^{-l(r)}$. Let $X \in \omega$ denote a discrete observable, a *statistical explanatory model* for $X$, is given by some unobservable, or uncontrolled variable $R \in [0, 1]$ with a probabilistic description given by a semimeasure $P_R$ over the unit interval $[0, 1]$ and some function $f$ such that $X = f(R)$. For some observation $x$ of the observable $X$, if $x = f(r)$, than $f, r$ is a probabilistic explanation of the observed data $x$, where $r$ represents the hidden or uncontrolled variables of the context where the value $x$ of $X$ is observed. The a-priori probability of occurence of $x$ is given by:

$$P_{f,R}(x) = \int dr\{r : x = f(r)\},$$

where Lebesgue integration over $r$, with respect to the measure $P_R$ is performed, and $f$ is assumed to be integrable.

For many contexts, it can be assumed that $f$ is partial computable and $P_R$ is enumerable. According to Lemma 2.1, without loss of generality, these assumptions are equivalent with assuming $P_R$ uniform over $[0, 1]$.

**Lemma 2.1.** *If the variable $R$ is distributed according to an enumerable $P_R$, and $f$ is partial computable, then there is a partial computable $f'$ and a uniform distributed variable $R'$ on $[0, 1]$, such that for all $x$:*

$$P_{f,R}(x) = P_{f',R'}(x).$$

*Proof.* First the function $\alpha$ is inductively defined. For any $x$, let $\alpha(x, 0) = 0$ and let

$$\alpha(x, t) = \alpha(2^n, t-1) + \sum\{P_t(z) - P_{t-1}(z) : z \leqslant x\}.$$

Remark that for every $r'$ such that

$$0 \leqslant r' \leqslant \sum\{P_{f,r}(x) : x \in 2^n\},$$

there is a unique $t$ such that $r' \in [\alpha(2^n, t-1), \alpha(2^n, t)]$. Therefore, each such $r'$ defines a unique $x$, such that

$$r' \in [\alpha(x-1, t), \alpha(x, t)].$$

If $l(r')$ is long enough, than also

$$[r'] \in [\alpha(x-1, t), \alpha(x, t)].$$

is satisfied, it is

$$r' \in [\alpha(x-1, t), \alpha(x, t) - 2^{l(r')}]. \tag{2}$$

Let $f'(r')$ be the function that is defined to be $x$ if there is an $x$ such that (2) is satisfied, and undefined otherwise. Remark that $f'$ is partial computable and satisfies the conditions of the Lemma. HOW CAN THIS BY NICELY WRITTEN OUT, HOW SHOULD I CHANGE DEFINITIONS ? $\square$

From now on, the variable $R$ will be assumed to have the uniform distribution, and $P_f$ is short for $P_{f,R}$. According to Proposition 2.2, the set of explanatory models is equivalent with the set of enumerable semimeasures.

**Proposition 2.2.** *For every partial enumerable $f$, the semimeasure $P_f$ is enumerable. For every enumerable semimeasure $P$, there is a partial computable function $f$ such that $P = P_f$.*

*Proof.* The first claim follows by definition. Let $\alpha(x, t)$ be as in the proof of Lemma 2.1. The second claim follows by choosing $f(r) = x$ if there is an $x$ such that for some $t$

$$\alpha(x, t) \leqslant r \leqslant \alpha(x, t) + 2^{-l(r)},$$

and $f(r) = \infty$ (undefined) otherwise. Remark that $f$ is partial computable and satisfies the conditions of the Lemma. $\square$

In [11, 12], the proof that high $K(K(x)|x)$ is rare, shows that for $t_i$ as defined there, which increases faster than any computable function of $i$, the probability for a prefix-free Turing machine that a program halts after time $t_i$ is bounded by $o(2^{-i})$. A similar proof shows that only for a small measure of hidden and uncontrolled variables $R$, there are $x$'es for which the exploratory model needs more computation time than $t_i$.

## 2.2 Causal explanations for a pair of observables

Different types of explanatory models are defined, and the corresponding universal elements are compared.

- An explanatory model for two discrete observables $X, Y$ is given by a partial computable function $f_{XY}$ and a variable $R$, uniformly distributed over $[0, 1]$, such that:

  $$(X, Y) = f_{XY}(R).$$

- An explanatory model for two *independent* discrete observables $X, Y$ is given by two partial computable functions $f_X, f_Y$ and two variables $R_X, R_Y$, independently and uniformly distributed over $[0, 1]$, such that:

$$\begin{aligned} X &= f_X(R_X) \\ Y &= f_Y(R_Y). \end{aligned}$$

- An explanatory model for two discrete variables $X, Y$ such that $X$ *causes* $Y$ is given by two partial computable functions $f_X, f_{Y|X}$ and two variables $R_X, R_Y$, independently and uniformly distributed over $[0, 1]$, such that

$$\begin{aligned} X &= f_X(R_X) & (3) \\ Y &= f_{Y|X}(X, R_Y). & (4) \end{aligned}$$

In a similar way as Proposition 2.2, the explanatory models defined above are equivalent with sets of semimeasures.

**Proposition 2.3.** *The universal elements of the sets of semimeasures corresponding to the explanatory models for $X, Y$, respectively independent $X, Y$, and $X$ causing $Y$, are given by $m(x, y)$, respectively, $m(x), m(y)$ and $m(x|y)m(y)$.*

*Proof.* This follows by the corresponding result of Proposition 2.2 and 1.2. □

Remark that the universal semimeasure $m(x, y)$ can be factorized.

**Lemma 2.4.** *Let $x^*$ be a program of length $K(x)$ that computes $x$, than*

$$m(x, y) =^* m(y|x^*)m(x).$$

*Proof.* Follows by applying the coding theorem and additivity of prefix-free Kolmogorov complexity [20]:

$$K(x) + K(y|x^*) =^+ K(x, y)$$

□

The corresponding test for the hypothesis that $x$ is independent from $y$ if $x$ is a probabilistic cause of $y$ is given by:

$$\frac{m(x)m(y|x)}{m(x)m(y)} = \frac{m(y|x)}{m(y)}.$$

The corresponding test for the hypothesis that $x$ is independent from $y$ if $x, y$ are generated in the most general way, is given by:

$$\frac{m(x, y)}{m(x)m(y)} = \frac{m(y|x^*)}{m(y)},$$

by Lemma 2.4. Remark that to approximate this test, a shortest model for $x$ might be needed. By Proposition 2.5 these tests can differ.

**Proposition 2.5.** *For every $n$ and all $x, y \in 2^n$*

$$\frac{m(y|x^*)}{m(y|x)} \leqslant^* n.$$

*For every n, there are x, y ∈ $2^n$ such that*

$$\frac{m(y|x^*)}{m(y|x)} \geqslant^* \frac{n}{\log n}$$

$$\frac{m(x|y)m(y)}{m(y|x)m(x)} \geqslant^* \frac{n}{\log n}$$

*Proof.* The first claim of Proposition 2.5 follows from the coding theorem [20]:

$$\log m(y|x^*) =^+ K(y|x^*)$$
$$\log m(y|x) =^+ K(y|x)$$

Remark by [1, Lemma 4.2] (see also appendix) that $K(x), x$ computes $x^*$ and by [20, page 242] it follows that

$$K(x^*|x) =^+ K(K(x)|x) \leqslant^+ \log n.$$

Remark that $K(y|x) \leqslant^+ K(y|x^*) + K(x^*|x)$. Combining the above equations shows the claim.

Now the second claim of Proposition 2.5 is shown. Remark that $K(x)$ can be computed from $x^*$, and that

$$K(x) =^+ K(K(x), x) =^+ K(K(x)) + K(x|K(x)^*).$$

Let $y = K(x)$. By applying the conditional coding theorem, it only needs to be shown that

$$K(K(x)|x) \geqslant^+ K(K(x)|x^*) + \log n - \log \log n$$

Remark that $K(K(x)|x^*) =^+ 0$. By [20, Theorem 3.8.1], it follows that for every $n$, there is at least one $x$ such that $K(K(x)|x) \geqslant^+ \log n - \log \log n$. □

**Corollary 2.6.** *There are hypothesis $S, T$ such that for every $n$, there are $x, y \in 2^n$ such that*

$$\frac{m^{(S \times T)^\uparrow}(x, y)}{m^{S^\uparrow \times T^\uparrow}(x, y)} \geqslant^* \frac{n}{\log n}$$

*Proof.* Let $S$ be the hypothesis that $x$ is generated by a partial computable function of a hidden variable $r_x$, and let $T$ be the hypothesis that $y$ is generated from $x$ by any function of a hidden variable $r_y$ and $x$. The universal element of $S^\uparrow$ is given by $m(x)$, the universal element of $T^\uparrow$ is given by $m(y|x)$. By Proposition 1.2, the universal element of $S^\uparrow \times T^\uparrow$ is given by $m(x)m(y|x)$.

It will now be shown that the universal element of $(S \times T)^\uparrow$ is given by $m(x, y)$. First remark that the semimeasure $m(x, y)/m(x)$ corresponds to some generalized structural equations where $y$ is generated from $x$ and a hidden variable $r_y$, by some function $f(x, r)$ that is not partial computable. Since $m(x, y)$ is also computable $m(x, y) \in (S \times T)^\uparrow$, but since $m(x, y)$ is also universal to the most general enumerable set of semimeasures, it must be universal to $(S \times T)^\uparrow$.

Finally it needs to be shown that for every $n$ there are $x, y \in 2^n$ such that:

$$\frac{m(x, y)}{m(x)m(y|x)}.$$

This follows from Proposition 2.5 and Lemma 2.4. □

Remark that if the partial computable functions $f_*$ in the definitions of the different exploratory models where chosen computable, the corresponding sets would in general not have a universal element. Any computable semimeasures $P(x|y), P(y)$, also generate semimeasures $P(y|x), P(x)$ that satisfy $P(x|y)P(y) = P(y|x)P(x)$. This means that it does not matter for the likelihood of $x, y$ to describe first $x$ and than $y$ and vice versa. This contrasts with the likelihood obtained from enumerable universal semimeasure by Proposition 2.5. Therefore, when $K(K(x)|x)$ is large, and $y$ contains much information about $K(x)$, the hypothesis that $x$ caused $y$ can be considered more plausible than vice versa.

Before a last type of hypothesis is introduced, the definition of total arguments of partial computable functions is given.

**Definition 2.7.** A partial computable function $f(x, y)$ on $2^{<\omega} \times U$ for some set $U$ has a *total argument* $x$, iff for any $x$ and $y \in U$ such that $f(x, y)$ is defined, also any $f(w, y)$ is defined with $w \in 2^{l(x)}$.

The hypothesis of $y$ being totally caused by $x$. is given by partial computable functions $f_X, f_{Y|X}$, such that $f_{Y|X}(X, R_X)$ has a total argument $X$, and two variables $R_X, R_Y$ independently and uniformly distributed over $[0, 1]$, satisfying equations (3) and (4). In a similar way, this hypothesis corresponds to a set of enumerable semimeasures, such that for each $n$, every $P(\epsilon|y)$ with $y \in 2^n$ is constant. These semimeasures have a universal element denoted as $m(y|\underline{x})$. Defining causality with this hypothesis can lead to fundamentally different results by the following Lemma.

**Lemma 2.8.** *For some $c$ and all $n \geqslant c$, there are $x, y \in 2^n$ such that:*

$$\log \frac{m(y|x)}{m(y|\underline{x})} \geqslant n - c \log n.$$

*Proof.* This follows from the standard coding theorem, the total coding theorem, and the existence of $x, y \in 2^n$ for any $n$ such that [2, 23]

$$K(x|y) - K(x|\underline{y}) \geqslant^+ n - 2 \log n.$$

For the definition of total conditional prefix-free complexity see [2, 23]. $\square$

This difference is due to Halting information present in $y$ [2]. It can be interpreted as follows: if the computation of $f_{Y|X}$ requires a time $t_i$ (see higher) for some large $i$, than $x, r$ must contain a large amount of Halting information. Remark that $t_i$ contains about $i$ bits of Halting information [3]. For a general partial computable function $f_{Y|X}$, this Halting information can be obtained from both arguments of the function $r$ and $x$. If $f_{Y|X}$ is total in its first argument, and $f_{Y|X}(x, r)$ is defined, than a program can be made that generates $t_i$ from $f_{Y|X}$ and $r$. Therefore if the computation of $y$ is so involved that it needs a time $t_i$, than $i$ bits of Halting information are present in $r$, and such probability decreases with $2^{-i}$. This is not the case for the partial computable $f_{Y|X}$.

## 2.3 Causal and influence-free explanations for two timeseries

Let $X, Y \in \omega^n$ be observables representing timeseries. The hypothesis that $X$ is an *instantaneous cause* of $Y$, is defined as the existence of partial computable functions $f_X, f_Y$, and variables $R_X, R_Y$ uniformly and independently distributed over $[0, 1]^n$ such that for all $i \leqslant n$:

$$\begin{aligned} X_i &= f_X(X^{i-1}, Y^{i-1}, R_X^i) \\ Y_i &= f_Y(X^i, Y^{i-1}, R_Y^i). \end{aligned}$$

See figure 1, right, black and red.

The hypothesis that $X, Y$ are *strict causal*, is defined as the existence of partial computable functions $f_X, f_Y$, and variables $R_X, R_Y$ uniformly and independently distributed over $[0, 1]^n$ such that for all $i \leqslant n$:

$$
\begin{aligned}
X_i &= f_X(X^{i-1}, Y^{i-1}, R_X^i) \\
Y_i &= f_X(X^{i-1}, Y^{i-1}, R_Y^i).
\end{aligned}
$$

See figure 2, black. Remark that by symmetry, if $X$ is a strict cause of $Y$, than $Y$ is a strict cause of $X$.

The hypothesis that $X$ is *influence-free* of $Y$, is defined as the existence of partial computable functions $f_X, f_Y$, and variables $R_X, R_Y$ uniformly and independently distributed over $[0, 1]^n$ such that for all $i \leqslant n$:

$$
\begin{aligned}
X_i &= f_X(X^{i-1}, R_X^i) \\
Y_i &= f_X(X^i, Y^{i-1}, R_Y^i).
\end{aligned}
$$

See figure 1, right, black and red.

The most general structure is obtained if hidden variables are shared. Therefore the hypothesis that $X, Y$ *can have hidden variables* is given by the partial computable functions $f_X, f_Y$ and the variable $R$ uniformly distributed over $[0, 1]^n$ such that for all $i \leqslant n$:

$$
\begin{aligned}
X_i &= f_X(R^i) \\
Y_i &= f_Y(R^i).
\end{aligned}
$$

This model is both equivalent with the models from figure 1, left, and 2, right, black and red.

## 2.4   Causal semimeasures

The hypothesis described in the previous subsection, correspond to sets of enumerable semimeasures which are investigated in this subsection. For $x \in 2^n$ and $i \leqslant n$, let

$$
P(x^i) = \sum \{ P(x^i v) : v \in 2^{n-i} \},
$$

and similar for $P(x^i, y^j)$ and $P(x^i | y)$. For $k \leqslant i \leqslant n$ and $l \leqslant j \leqslant n$, let

$$
P(x^i, y^j | x^k, y^l) = \frac{P(x^i, y^j)}{P(x^k, y^l)}.
$$

**Definition 2.9.** Let $x, y \in 2^n$.

- The causal semimeasure and the instantaneous causal semimeasure, *associated* with a conditional semimeasure $P(x|y)$ is given by:

$$
\begin{aligned}
P(x|y \uparrow) &= \prod \{ P(x_i | x^{i-1}, y^{i-1}) : i \leqslant n \} \\
P(x|y \uparrow^+) &= \prod \{ P(x_i | x^{i-1}, y^i) : i \leqslant n \}.
\end{aligned}
$$

- A conditional semimeasure $P(y|x)$ is *causal* respectively *instantaneous causal*, iff for all $i \leqslant n$ respectively

$$
\begin{aligned}
P(y|x \uparrow) &= P(y|x) \\
P(y|x \uparrow^+) &= P(y|x).
\end{aligned}
$$

- $x$ is *influence-free* of $y$ according to a semimeasure $P(x, y)$, iff:

$$P(x|y \uparrow) = P(x), \tag{5}$$

when defined.

**Proposition 2.10.** *For any semimeasure $P(x, y)$, the following statements are equivalent:*

(i) $P(y|x)$ *is instantaneous causal.*

(ii) $\forall i \leqslant n \forall x, y \in 2^n \big[ P(y^i|x) = P(y^i|x^i) \big]$ *where defined.*

(iii) $\forall i \leqslant n \forall x, y \in 2^n \big[ P(x|x^i, y) = P(x|x^i) \big]$ *where defined.*

(iv) $\forall i \leqslant n \forall x, y \in 2^n \big[ P(x_{i+1}|x^i, y) = P(x_{i+1}|x^i) \big]$ *where defined.*

(v) $x$ *is influence-free of $y$ according to $P(x, y)$.*

*Proof.* $(i) \rightarrow (ii)$:
Let

$$P(y^i|x^i \uparrow^+) = \prod \{ P(y_j|y^{j-1}, x^j) : j \leqslant i \}.$$

First it is shown that

$$P(y^i|x) = P(y^i|x^i \uparrow^+). \tag{6}$$

Suppose that for some $y^i, x$: $P(y^i|x) > P(y^i|x^i \uparrow^+)$, than for every $j = i + 1, ..., n$, choose $y_j$ such that $P(y_{j+1}|x, y^j) \geqslant P(y_{j+1}|x^{j+1}, y^j)$. Remark that this is always possible. This shows that:

$$
\begin{aligned}
P(y|x) &= P(y^i|x) \prod \{ P(y_j|x, y^{j-1}) : j = i + 1 ... n \} \\
&> P(y^i|x^i \uparrow) \prod \{ P(y_j|x^j, y^{j-1}) \} \\
&= P(y|x \uparrow^+).
\end{aligned}
$$

which contradicts $(i)$. $(ii)$ follows by

$$
\begin{aligned}
P(y^i|x^i) &= \sum_x P(y^i|x) P(x|x^i) \\
&= \sum_x P(y^i|x^i \uparrow^+) P(x|x^i) \\
&= P(y^i|x^i \uparrow^+).
\end{aligned}
$$

$(ii) \rightarrow (iii)$: By Bayes theorem.
$(iii) \rightarrow (iv)$: By summing over all $x^{i+1}v$ with $v \in 2^{n-i-1}$.
$(iv) \rightarrow (v)$: By definition.
$(v) \rightarrow (i)$: By remarking that

$$P(x|y \uparrow) P(y|x \uparrow^+) = P(x, y) = P(x) P(y|x).$$

$\square$

Remark that the set of causal, and instantaneous causal semimeasures is testable and convex. Therefore, they have universal elements $m(x|y \uparrow), m(x|y \uparrow^+)$.

The introduced hypotheses from the previous subsection correspond to sets of enumerable semimeasures which have a universal element. This follows by the same argument as Proposition 2.5. The corresponding universal elements are given by Proposition 2.11.

**Proposition 2.11.** *The universal element of the hypothesis that*

1. *$X$ is an instantaneous cause of $Y$, is given by $m(x|y \uparrow)m(y|x \uparrow^+)$.*

2. *$X, Y$ are strict causal, is given by $m(x|y \uparrow)m(y|x \uparrow)$.*

3. *$X$ is influence-free of $Y$, is given by $m(x)m(y|x \uparrow)$.*

4. *$X, Y$ have hidden common variables, is given by $m(x, y)$.*

*Proof.* The corresponding sets of semimeasures are products of convex enumerable sets of enumerable semimeasures. The result follows by Proposition 1.2. $\square$

The universal elements define ideal hypotheses tests. Some of them can be simplified within a constant factor, using $m(x, y) =^* m(x)m(y|x^*)$.

- Suppose that $X$ is an instantaneous cause of $Y$,
  are $X, Y$ strict causal according to data $x, y$?
  Figure 2, left.
  $$\frac{m(x|y \uparrow)m(y|x \uparrow^+)}{m(x|y \uparrow)m(y|x \uparrow)} = \frac{m(y|x \uparrow^+)}{m(y|x \uparrow)} \tag{7}$$

- Suppose that $X, Y$ are strict causal,
  is $Y$ influence-free of $X$ according to data $x, y$?
  Figure 1, right.
  $$\frac{m(x|y \uparrow^+)m(y|x \uparrow)}{m(x|y \uparrow^+)m(y)} = \frac{m(y|x \uparrow)}{m(y)} \tag{8}$$

- Suppose $X, Y$ can have hidden variables,
  is $X$ an instantaneous cause of $Y$ according to data $x, y$ ?
  $$\frac{m(x, y)}{m(x|y \uparrow^+)m(y|x \uparrow)} \tag{9}$$

- Suppose $X, Y$ can have hidden variables,
  are $X, Y$ strict causal ?
  Figure 2, right.
  $$\frac{m(x, y)}{m(x|y \uparrow)m(y|x \uparrow)} \tag{10}$$

- Suppose $X, Y$ can have hidden variables,
  is $Y$ influence-free of $X$ ?
  $$\frac{m(x|y^*)m(y)}{m(x|y \uparrow^+)m(y)} = \frac{m(x|y^*)}{m(x|y \uparrow^+)} \tag{11}$$

The significances of ideal independence tests can now be decomposed as the product of the significances of the tests above. For example as the product of the tests of equations (10), (8), and (8) applied to $x$, or as the decomposition. (9), (7), (8), and (8) applied to $x$, or as the decomposition.

$$\begin{aligned}
\frac{m(x, y)}{m(x)m(y)} &= \frac{m(x, y)}{m(x|y \uparrow)m(y|x \uparrow)} \frac{m(y|x \uparrow)}{m(y)} \frac{m(x|y \uparrow)}{m(x)} \\
&= \frac{m(x, y)}{m(x|y \uparrow^+)m(y|x \uparrow)} \frac{m(y|x \uparrow^+)}{m(y|x \uparrow)} \frac{m(y|x \uparrow)}{m(y)} \frac{m(x|y \uparrow)}{m(x)}
\end{aligned} \tag{12}$$

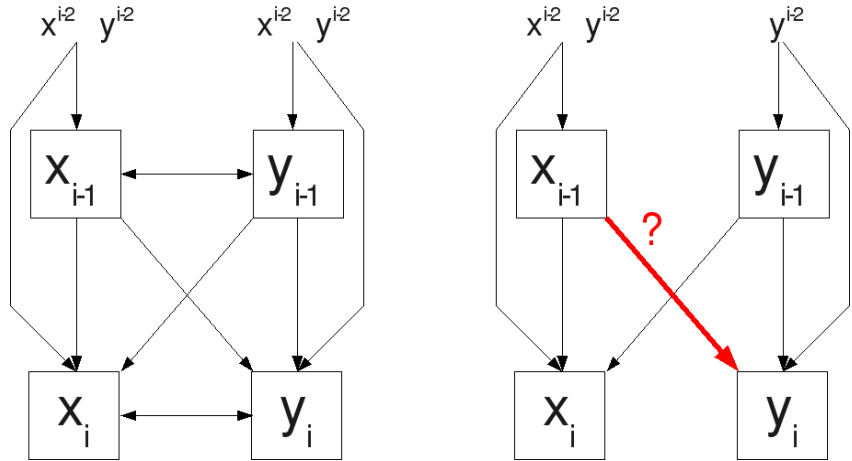

Figure 1: Left: general system. Right: suppose that $X, Y$ are strict causal, is $Y$ influence-free of $X$ ?

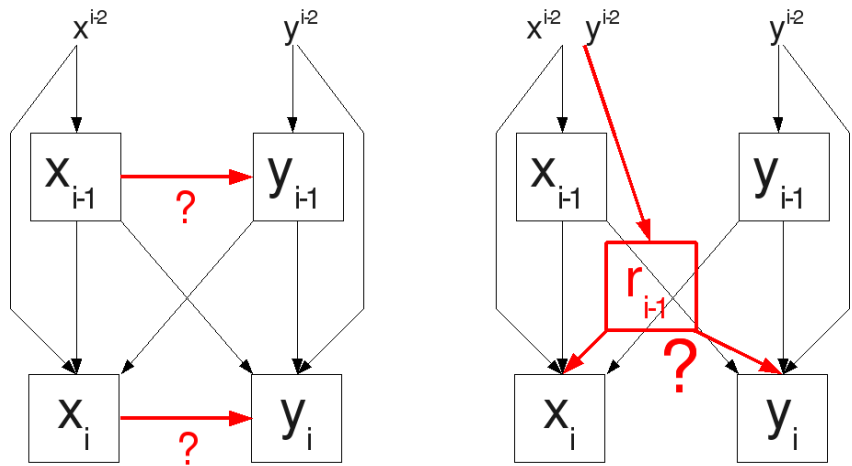

Figure 2: Left: Suppose that $X$ is an instantaneous cause of $Y$, are $x, Y$ strict causal ? Right: Suppose $X, Y$ can have hidden variables, are $X, Y$ strict causal ?

## 2.5 Information transfer and instantaneous information transfer

Equation (12) allows a nice information theoretic interpretation. Let

$$
\begin{array}{rcl}
I(x;y) & = & \log \dfrac{m(x,y)}{m(x)m(y)} \\[2ex]
EIT(x \leftarrow y) & = & \log \dfrac{m(x \lceil y \uparrow)}{m(x)} \\[2ex]
EIT(x \uparrow; y \uparrow) & = & \log \dfrac{m(x,y)}{m(x|y \uparrow)m(y|x \uparrow)}.
\end{array}
$$

Equation (12) becomes now:

$$
I(x;y) = EIT(x \leftarrow y) + EIT(y \leftarrow x) + EIT(x \uparrow; y \uparrow).
$$

Suppose that $x, y$ have no instantaneous connections, than the mutual information of $x, y$, can be considered as the sum of information flowing from the past of $x$ to $y$, from the past of $y$ to $x$, and information obtained by $x, y$ through a hidden source.

However, a decomposition of mutual information as the sum of information flowing from the past and the present of $x$ to $y$ and from the past of $y$ to $x$ is not possible. Also in this setting, there can be a different instantaneous information flow if it is assumed that information is instantaneously flowing from $x$ to $y$, or from $y$ to $x$. Both claims follow from Proposition 2.12.

**Proposition 2.12.** *For every $n$ there exist $x, y \in 2^n$ such that:*

$$
EIT(x \uparrow; y \uparrow) - \log \frac{m(y|x \uparrow^+)}{m(y|x \uparrow)} \geqslant o(n).
$$

*For there exist $x, y \in \omega^n$ such that:*

$$
\log \frac{m(x|y \uparrow^+)}{m(x|y \uparrow)} - \log \frac{m(y|x \uparrow^+)}{m(y|x \uparrow)} \geqslant o(\sum \{\log x_i + \log y_i : i \leqslant n\}).
$$

*Proof.* The first claim follows from Lemma 2.13, the second claim follows by an analogue result of Lemma 2.13 and from [1, Proposition ...][2]     □

**Lemma 2.13.** *There is a constant $c$ such that for any $n \geqslant \frac{1}{c}$, there exist $x, y \in 2^n$ such that*

$$
\log \frac{m(x,y)}{m(x|y \uparrow)m(y|x \uparrow^+)} \geqslant cn.
$$

*Proof.* In [6], on-line decision complexity $KR$ and on-line a priori complexity $KA$ are defined and it is shown that for the task $T = \epsilon \rightarrow x_1; y_1 \rightarrow x_2; ...; y_{n-1} \rightarrow x_n$:

$$
KR(T) \leqslant^+ KA(T) + 2 \log KR(T). \tag{13}
$$

Remark that in this result all quantities can be conditioned on there length, and that $KA$ can be replaced by $-\log m(x|y \uparrow)$ for any enumerable on-line semimeasure. The analogue for $-\log m(x|y \uparrow^+)$ also holds. Let $K(x|y \uparrow)$ and $K(x|y \uparrow^+)$ be length conditional on-line decision complexities as defined in [1]. On-line decision complexity

---

[2]The current draft version of the paper from August 2009 contains only the proof for the result in Lemma 2.13, however, the proof can be reformulated such that also this result is shown.

and length conditional on-line decision complexity are related within $2\log n$ terms. Therefore:

$$
\begin{aligned}
K(x|y\uparrow) &\leqslant -\log m(x|y\uparrow) + O(\log n) \\
K(y|x\uparrow^+) &\leqslant -\log m(y|x\uparrow^+) + O(\log n).
\end{aligned}
$$

By the coding theorem:

$$K(x,y) =^+ -\log m(x,y) \tag{14}$$

By the main result of [1], there exist for each $n$, $x,y \in 2^n$ such that:

$$K(x|y\uparrow) + K(y|x\uparrow^+) - K(x,y) \geqslant o(n).$$

Combining the above equations shows that Lemma. □

## 3 Associated causal semimeasures

In the previous section, enumerable causal semimeasures where derived as corresponding to structural equations with partial computable functions. In this section causal semimeasures are investigated that are in a Bayesian way associated to enumerable conditional semimeasures and enumerable bivariate semimeasures.

**Definition 3.1.** A semimeasure $P(x|y)$ is *associated causal* respectively *associated instantaneous causal* if it there is an enumerable conditional semimeasure $Q(x|y)$ such that $P(x|y) = Q(x|y\uparrow)$ respectively $P(x|y) = Q(x|y\uparrow^+)$.

Remark that an enumerable causal semimeasure $P(x|y)$ is associated causal, since it equals its own association. Also remark that the set of associated causal semimeasures is not convex. Since with any bivariate semimeasure $P(x,y)$, a conditional semimeasure is associated, one can associate a causal and instantaneous causal semimeasures also with $P(x,y)$.

**Lemma 3.2.**

$$P(x,y) = P(\epsilon,\epsilon)P(x|y\uparrow)P(y|x\uparrow^+). \tag{15}$$

*Proof.* Remark that for the causal semimeasures $P(x|y\uparrow)$ and $P(y|x\uparrow^+)$ associated with $P(x,y)$ one has

$$
\begin{aligned}
P(x|y\uparrow) &= \frac{P(x^1,y^0)}{P(x^0,y^0)} \cdots \frac{P(x^n,y^{n-1})}{P(x^{n-1},y^{n-1})} \\
P(x|y\uparrow^+) &= \frac{P(x^1,y^1)}{P(x^0,y^1)} \cdots \frac{P(x^n,y^n)}{P(x^{n-1},y^n)}.
\end{aligned}
$$

□

### 3.1 Non existence of universal elements

**Proposition 3.3.** *If $P(x|y\uparrow)$ is a causal semimeasure associated with an enumerable semimeasure $P(x|y) > 0$, than there exists an enumerable semimeasure $Q(x|y)$ and $x,y \in 2^n$ such that*

$$\log \frac{Q(x|y\uparrow)}{P(x|y\uparrow)} > o(n). \tag{16}$$

*Proof of Proposition 3.3 first part: definition of Algorithm 1.*

Let $N = 2n$ and let the set $2^n \times 2^n$ be associated with $2^N$ by mapping $x, y$ to $z = x_1 y_1 x_2 y_2 \ldots x_n y_n$. Let with abuse of notation $P$ denote the restriction of $P$ on $2^N$. For any restricted semimeasure $P$ and $v \in 2^i, i \leqslant N$, let $P(v\ldots)$ denote the restriction of $P$ on the strings $vu$ for all $u \in 2^{N-i}$. For $b \in \{0, 1\}$, let $\bar{b} = 1 - b$. The strings of $2^N$ can be considered as branches in a tree. For $z \in 2^N$, $z$ is a local minimal branch, iff it satisfies for all $i \leqslant n$:

$$P(z^i) \leqslant P(z^{i-1}\overline{z_i}).$$

For a local minimal branch $z$, the nodes $z^{2i+1}\overline{z^{2i+2}}$ for $i \leqslant n-1$ are called load nodes. Algorithm 1 generates for every restriction $P$ on $2^n$ of a computable semimeasure a computable semimeasure $Q$ on $2^n$ such that all leafs $w$ have half weight, it is $Q(w) = P(w)/2$, except for leafs descending from load nodes which have $Q(w) = P(w)$. This implies that the weights of the uneven local minimal nodes are proportionally more heavy than the weights of the even local minimal nodes, which shows the result of Lemma 3.4.

**Lemma 3.4.** *If $P(x|y \uparrow)$ is a causal semimeasure associated with a computable semimeasure $P(x, y) > 0$, than $Q = grow(P)$, with algorithm grow defined in Algorithm 1, is computable and satisfies:*

$$\log \frac{Q(x|y \uparrow)}{P(x|y \uparrow)} \;>\; o(n).$$

*Proof.* Algorithm 1 constructs $Q$ from $P$ such that:

$$Q(w) = \begin{cases} P(w) & \text{if } w \text{ is a load leaf,} \\ \frac{1}{2}P(w) & \text{otherwise.} \end{cases}$$

For $i < N$ and $z$ the local minimal leaf,

$$P(z^{i+1}) \leqslant \frac{1}{2}P(z^i),$$

and for $i < n$,

$$
\begin{aligned}
Q(z^{2i}) - \frac{1}{2}P(z^{2i}) &= Q(z^{2i+1}) - \frac{1}{2}P(z^{2i+1}) \\
&\geqslant \frac{1}{2}P(z^{2i+1}\overline{z_{2i+2}}) \geqslant \frac{1}{4}P(z^{2i+1}).
\end{aligned}
$$

This shows that:

$$
\begin{aligned}
\frac{Q(z^{2i+1})}{Q(z^{2i})} &= \frac{\frac{1}{2}P(z^{2i+1}) + Q(z^{2i+1}) - \frac{1}{2}P(z^{2i+1})}{\frac{1}{2}P(z^{2i}) + Q(z^{2i}) - \frac{1}{2}P(z^{2i})} \\
&= \frac{\frac{1}{2}P(z^{2i+1}) + Q(z^{2i+1}) - \frac{1}{2}P(z^{2i+1})}{\frac{1}{2}P(z^{2i}) + Q(z^{2i+1}) - \frac{1}{2}P(z^{2i+1})} \\
&\geqslant \frac{\frac{1}{2}P(z^{2i+1}) + \frac{1}{4}P(z^{2i+1})}{\frac{1}{2}P(z^{2i}) + \frac{1}{4}P(z^{2i+1})} \\
&\geqslant \frac{P(z^{2i+1})}{P(z^{2i})} \frac{1 + \frac{1}{2}}{1 + \frac{1}{2}P(z^{2i+1})/P(z^{2i})} \\
&\geqslant \frac{6}{5} \frac{P(z^{2i+1})}{P(z^{2i})}.
\end{aligned}
$$

$$
\begin{aligned}
-\log Q(x|y \uparrow) &= \sum_{i \leqslant n} -\log \frac{Q(z^{2i+1})}{Q(z^{2i})} \\
&= \sum_{i \leqslant n} -\log \frac{P(z^{2i+1})}{P(z^{2i})} - \log \frac{6}{5} \\
&\geqslant -\log P(x|y \uparrow) - n \log \frac{6}{5}.
\end{aligned}
$$

Remark that Algorithm 1 constructs a computable $Q$ from a computable $P$. $\qquad\square$

**Data**: $P$
**Result**: $Q$
**begin**
$\quad z \longleftarrow$ local minimal branch in $2^N$
$\quad Q \longleftarrow \frac{1}{2}P$
$\quad$ **for** $i$ *from* $0$ *to* $n-1$ **do**
$\quad\quad Q(z^{2i+1}\overline{z^{i+2}}...) \longleftarrow P(z^{i+1}\overline{z^{i+2}}...)$ $\qquad$ (load node)
**end**

**Algorithm 1**: grow

**Data**: $P_t$
**Result**: $Q_t$
**begin**
$\quad \nu = \min\{P_0(w) : w \in 2^N\}$
$\quad Q_0(w) \longleftarrow$ grow$2^{-1/\nu}P_0(w)$
$\quad S \longleftarrow P_0(\epsilon)$
$\quad s \longleftarrow 0$
$\quad$ **for** $t$ *from* $0$ *to* $\infty$ **do**
$\quad\quad$ **if** $P_t(\epsilon) - S > \nu$ **then** stage $s$: new $Q_t$ is grown
$\quad\quad\quad S \longleftarrow P_t(\epsilon)$
$\quad\quad\quad s \longleftarrow s + 1$
$\quad\quad\quad Q_t \longleftarrow$ grow$(2^{-1/\nu+s}P_t)$
$\quad\quad$ **else**
$\quad\quad\quad Q_t \longleftarrow \frac{P_t}{P_{t-1}}Q_{t-1}$
**end**

**Algorithm 2**: grow_semimeasure

*Proof of Proposition 3.3 second (last) part.*
The causal semimeasure associated with a bivariate semimeasure $P(x, y)$, is the causal semimeasure associated with the conditional semimeasure $P(x|y)$. For every enumerable conditional semimeasure $P(x|y) > 0$, there is an enumerable semimeasure $Q(x, y) > 0$ such that $Q(x|y) = P(x|y)$. Therefore, to show Proposition 3.3, it suffices to show the proposition for causal semimeasures associated with bivariate semimeasures.

Let $\nu = \min\{P_0(w) : w \in 2^N\}$. Remark that the enumeration $P_t$ can be chosen such that $\nu > 0$, since $P(x, y) > 0$. Algorithm 2 uses Algorithm 1, to define an enumerable semimeasure $Q_t$ from the enumerable semimeasure $P_t$. It is now shown that $Q_{t+1} \geqslant Q_t$: suppose that $t, t + 1$ are in the same stage $s$, than this is easily observed, if at time $t + 1$ a new stage $s + 1$ is reached, a new $Q_{t+1}$ is grown from $P_{t+1}$ multiplied with a factor $2^{1/\nu+s+1}$, which is doubled relative to the previous stage.

Therefore if $w$ was a non-load leaf at time $t$, and becomes a load leaf at time $t+1$, one still has $Q_{t+1}(w) \geqslant Q_t(w)$. By Lemma 3.4 it follows for every $t$ that initiates a new stage, that:

$$\frac{Q_t(x|y\uparrow)}{P_t(x|y\uparrow)} \geqslant o(n).$$

By Lemma 3.5, this equation also hold for the $t$ subsequent to the $t$'s initiating a new stage. Therefore, the equation holds for any $t$. $\qquad\square$

**Lemma 3.5.** *Suppose that for some $\nu > 0$, and some semimeasures $P, Q \in 2^N$, satisfy for all $w \in 2^N$,*

$$\begin{aligned}
P(w) &\geqslant \nu \\
Q(w) &\geqslant P(w) \\
Q(\epsilon) &\leqslant P(\epsilon) + \nu,
\end{aligned}$$

*than*

$$\frac{1}{2} \leqslant \frac{Q(x|y\uparrow)}{P(x|y\uparrow)} \leqslant 2.$$

*Proof.* Remark that for any $i \leqslant N$ and $w \in 2^i$, one has $Q(w) \leqslant P(w) + \nu$ and $P(w) \geqslant \nu 2^{N-i}$. Since any branch of depth $j$ has $2^{N-j}$ leafs, one has that $P(w^j) \geqslant \nu 2^{N-j}$.

$$\begin{aligned}
Q(x|y\uparrow) &= \prod\{\frac{Q(w^{2i+1})}{Q(w^{2i})} : i < n\} \\
&\geqslant \prod\{\frac{P(w^{2i+1})}{P(w^{2i}) + \nu} : i < n\} \\
&\geqslant \prod\{\frac{P(w^{2i+1})}{P(w^{2i}) + P(w^{2i})2^{-N+2i}} : i < n\} \\
&\geqslant P(x|y\uparrow)\prod\{\frac{1}{1 + 2^{-2i}} : i < n\} \\
&\geqslant \frac{1}{2}P(x|y\uparrow).
\end{aligned}$$

$$\begin{aligned}
Q(x|y\uparrow) &= \prod\{\frac{Q(w^{2i+1})}{Q(w^{2i})} : i < n\} \\
&\leqslant \prod\{\frac{P(w^{2i+1}) + \nu}{P(w^{2i})} : i < n\} \\
&\geqslant \prod\{\frac{P(w^{2i+1}) + P(w^{2i+1})2^{-N+2i+1}}{P(w^{2i})} : i < n\} \\
&\geqslant P(x|y\uparrow)\prod\{1 + 2^{-2i-1} : i < n\} \\
&\geqslant 2P(x|y\uparrow).
\end{aligned}$$

$\qquad\square$

## 3.2 Causal semimeasures associated with a universal semimeasure

Let $m(x\lceil y\uparrow)$ be the causal semimeasure associated with $m(x, y)$ and let $\tilde{m}(x\lceil y\uparrow)$ be the causal semimeasure associated with $m(x|y)$.

**Proposition 3.6.** *There is a constant $c > 0$ such that for all $n \geqslant \frac{1}{c}$, there are $x, y \in 2^n$ with*

$$\log \frac{m(x \lceil y \uparrow)}{m(x|y \uparrow)} \quad \geqslant \quad cn,$$

*and there are $x, y \in 2^n$ with*

$$\log \frac{\tilde{m}(x \lceil y \uparrow)}{m(x|y \uparrow)} \quad \geqslant \quad cn.$$

*Proof.* The first claim of Proposition 3.6 is now proved. Let $m_t(x \lceil y \uparrow)$ be associated with $m_t(x, y)$. Suppose the first claim is not true, than for any constant $c > 0$ and any $x, y \in 2^n$:

$$m(x \lceil y \uparrow) \leqslant 2^{cn} m(x|y \uparrow).$$

Fix some enumeration $m_t(x|y \uparrow)$ of $m(x|y \uparrow)$. It can now be assumed that for some enumeration $m_t(x, y)$ of $m(x, y)$, one has:

$$m_t(x \lceil y \uparrow) \leqslant 2^{cn} m_t(x|y \uparrow).$$

Let

$$
\begin{aligned}
\mu_t(x|y) &= 2^{cn} m_t(x|y \uparrow) \\
\mu_t(y|x \uparrow^+) &= \frac{m_t(x, y)}{\mu_t(x|y \uparrow)}.
\end{aligned}
\tag{17}
$$

Remark that $\mu_t(x|y)$ may not be a semimeasure, but $\mu_t(y|x \uparrow^+)$ is a causal semimeasure by Lemma 3.2 and using $m(\epsilon, \epsilon) \leqslant \frac{1}{O(1)}$. The $\mu$ functions are limit-computable, and are used to define the following $\nu$ functions:

$$
\begin{aligned}
s_t &= \mathrm{argmax}_s \{ \mu_s(y|x \uparrow^+) : s \leqslant t \} \tag{18} \\
\nu_t(y|x \uparrow^+) &= \mu_{s_t}(y|x \uparrow^+) \tag{19} \\
\nu_t(x|y \uparrow) &= \mu_{s_t}(x|y \uparrow) \frac{m_t(x, y)}{m_{s_t}(x, y)}. \tag{20}
\end{aligned}
$$

Also, remark that $\nu_t(y|x \uparrow^+)$ is increasing in $t$. $\nu_t(x|y \uparrow)$ is also increasing in $t$ since Equation (20) shows that if $s_t = s_{t+1}$ than $\nu_t(x|y \uparrow) \leqslant \nu_{t+1}(x|y \uparrow)$, and for any $t$ such that $s_t < s_{t+1}$, remark that $s_{t+1} = t + 1$ and by equation (20):

$$\nu_t(x \lceil y \uparrow) \leqslant \mu_t(x|y \uparrow) \leqslant \mu_{t+1}(x|y \uparrow) = \nu_{t+1}(x \lceil y \uparrow).$$

Remark that using equations (18) and (17) show that:

$$
\begin{aligned}
\nu_t(x|y \uparrow) &= \mu_{s_t}(x|y \uparrow) \frac{m_t(x, y)}{m_{s_t}(x, y)} \\
&= \mu_t(x|y \uparrow) \frac{m_t(x, y)}{\mu_t(x, y)} \frac{\mu_{s_t}(y|x \uparrow^+)}{m_{s_t}(y|x \uparrow^+)} \\
&\leqslant \mu_t(x|y \uparrow).
\end{aligned}
$$

This shows that:

$$
\begin{aligned}
m(x, y) &= \nu(x|y \uparrow) \nu(y|x \uparrow^+) \\
&\leqslant 2^{nc} m(x|y \uparrow) \nu(y|x \uparrow^+) \\
&\leqslant^* 2^{nc} m(x|y \uparrow) m(y|x \uparrow^+).
\end{aligned}
$$

This equation is valid for any $c$ and $x, y \in 2^n$, contradicting Lemma 2.13. This shows the first claim of Proposition 3.6.

The second claim of the Proposition follows by replacing $m_t(x, y)$ by

$$\tilde{m}(x, y) = m(x|y)m(y).$$

Remark that $\tilde{m}(x \lceil y \uparrow)$ is also the causal associated semimeasure of $\tilde{m}(x, y)$. The analogue contradiction of Equation (3.2) is derived by remarking that by Proposition 2.5:

$$m(x, y) \leqslant^* m(x|y)m(y)n.$$

$\square$

**Corollary 3.7.** *$m(x \lceil y \uparrow)$ and $\tilde{m}(x \lceil y \uparrow)$ are not enumerable.*

**Question 3.8.** *Let $S$ be the set of causal semimeasures associated with an universal enumerable semimeasure. How much can two elements of $S$ differ. Has $S$ a universal element ?*

The relations betweens the sets of associated and enumerable causal semimeasures are represented in figure 3.

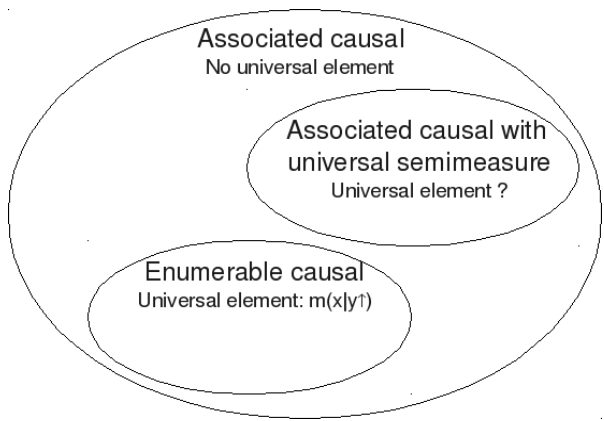

Figure 3: Relations between sets of causal semimeasures and existence of universal elements.

## 3.3 Associated information transfer and instantaneous common information

For any universal semimeasure $m$, the associated causal semimeasure will be denoted as $m(x \lceil y \uparrow)$. Associated information transfer and instantaneous common information are given by:

$$
\begin{aligned}
AIT(x \leftarrow y) &= -\log \frac{m(x \lceil x \uparrow)}{m(x)} \\
AIT(x \uparrow; y \uparrow) &= -\log \frac{m(x, y)}{m(x|y \uparrow)m(y|x \uparrow)}.
\end{aligned}
$$

Remark that:

$$I(x; y) = AIT(x \leftarrow y) + AIT(y \leftarrow x) + AIT(x \uparrow; y \uparrow).$$

Associated simultaneous information transfer has also another representation.

$$AIT(x \uparrow; y \uparrow) = \log \frac{m(x \lceil y \uparrow^+)}{m(x \lceil y \uparrow)} = \log \frac{m(y \lceil x \uparrow^+)}{m(y \lceil x \uparrow)}.$$

This means that in contrast with the enumerable instantaneous common information, the associated instantaneous common information can be either interpreted as an instantaneous information flow flowing from $x$ to $y$, or flowing from $y$ to $x$, or simultaneously flowing from a hidden source to $x$ and $y$.

# 4 Shannon information transfer and minimal sufficient statistics

## 4.1 Granger causality and Shannon information transfer

Statistical tests used in engineering literature can often be structured as follows: first a model is fitted on the data and than influence is derived from:

- Some parameters in the model, as for example by the use of directed transfer functions [16] and partial directed coherences [27].

- The complexity or magnitude of the noise of the data relative to the model, as for example with the use of Granger Causality [15, 10, 14, 8, 16], and Shannon information transfer [24, 26, 28, 32].

By an on-line coding theorem [6], the ideal statistical tests based on enumerable information transfer, can be informally assumed to derive influence from the sum of the complexity of the model, and the complexity of the noise. It is not clear whether such algorithms perform better [5].

Let $E(X^+|X^-)$ denote some average error of a prediction strategy of observations of the observable $X$ given its past observations. Let $E(X^+|X^-, Y^-)$ be similar where the prediction strategy also uses the past of $Y$. In its most general form [10, 14], $Y$ is said to Granger causal $X$ iff

$$E(X^+|X^-) - E(X^+|X^-Y^-))$$

is large. The most common choice for $E(.|.)$ to define Granger causality is the mean squared error.

Another choice for $E(.|.)$ is Shannon entropy. The following expressions provide definitions for Shannon mutual information, information transfer and instantaneous mutual information:

$$
\begin{aligned}
SI_P(X; Y) &= \sum_{x,y} P(x, y) \log \frac{P(x, y)}{P(x) P(y)} \\
SIT_P(X \leftarrow Y) &= \sum_{x,y} P(x, y) \log \frac{P(x|y \uparrow)}{P(x)} \\
SIT_P(X \uparrow; Y \uparrow) &= \sum_{x,y} P(x, y) \log \frac{P(x, y)}{P(x|y \uparrow) P(y|x \uparrow)}.
\end{aligned}
$$

Remark that

$$SI_P(x; y) = SIT_P(x \leftarrow y) + SIT_P(y \leftarrow x) + SIT_P(x \uparrow; y \uparrow).$$

A general procedure of deriving influence for the procedures in [24, 26, 28, 32] is given by fitting some models $P(X_t|X_{t-k...t-1})$, $P(X_t|X_{t-k...t-1}, Y_{t-k...t-1})$, to the corresponding data segments and similar for $Y_t$, and finally computing the statistic $SIT_P(X \leftarrow Y) - SIT_P(X \leftarrow Y)$. A confidence for the sign of the statistic can be obtained by running the procedure on some randomized permutation of the sequences $x$ and $y$.

The continuous entropy of a Normal distribution is given by $\sqrt{2\pi e}\sigma$ [29]. This implies that when the error of the observed data relative to some model is assumed to be Normal distributed, the Shannon entropy is estimated by the root mean squared error, in correspondence with common definitions of Granger causality.

When it is assumed that $P$ is a good model for the data, in a frequentistic interpretation, this means that for repetitive observation of the data, the data is distributed according to $P$, then ideal on-line data compression is with overwhelming probability performed by Shannon-Fano code [20]. The expected difference of the code-length of the on-line Shannon-Fano code and the unconditioned code is given by the Shannon information transfer. The expected code-length for optimal on-line encoding is given by enumerable information transfer within small terms [2]. Therefore, mean Shannon information transfer and mean enumerable information transfer are equal within some constant. A formal version of this statement is given by Proposition 4.1.

**Proposition 4.1.**

$$\sum P(x,y)EIT(x \leftarrow y) \quad =^+ \quad \sum P(x,y)AIT(x \leftarrow y) \quad =^+ \quad SI_P(X \leftarrow Y) \pm K(P)$$

*Proof.* The right equality has the same structure as the proof of [13, Lemma II.4]. The left equality is shown in [2]. □

## 4.2 Minimal sufficient statistics and ideal Shannon information transfer

Algorithms for extracting $P$ from $x, y$ as in the previous subsection, are often designed to let $P$ model as much as possible properties that appear frequently within the signal, while at the same time keeping the descriptional complexity of $P$ low.

To idealize this procedure, it has been conjectured [18] that in the case of multivariate models a, the constructed $P$ should be chosen as a probabilistic minimal sufficient statistic of the data $x, y$ [13]. Two ways of assigning causal relationships from such a $P$ exists, either by computing $SIT_P$ or by extracting a graphical schema from $P$. In [18] it is argued informally that for the multivariate case, a minimal Bayesian network is a minimal sufficient statistic. At [18, Lemma 4], it is claimed that if a two-part code satisfying some syntactical form results in an incompressible string, the first part is the probabilistic minimal sufficient statistic. However, it is argued here that in many cases, any plausible graphical causal representation cannot be contained in a probabilistic minimal sufficient statistic of the data. If the plain Kolmogorov complexity of the weak minimal sufficient statistic [3] is computable from $n$, for example, it is relatively small, and the causal structure of the observables is not equivalent with some initial segment of the Halting sequence, than the causal information can not be in the minimal sufficient statistic. To show this, assume that the plain complexity of the weak minimal sufficient statistic is computable from $n$, it follows from $C(P) =^+ K(P|C(P)) =^+ K(P)$ [20, Lemma 3.1] that the weak minimal sufficient statistic is also a minimal sufficient statistic. From [3, Proposition 6.5 and Proposition 7.6] it follows that every weak minimal sufficient statistic is equivalent with some initial segment of the Halting sequence. Therefore, the minimal sufficient statistic can not include the causal information.

In the proof of Proposition 4.2, a more formal argument will be given for the bivariate case, by providing a more explicit construction of two pairs of strings where a very different plausible causal relationship is present but also have the same minimal sufficient statistic. For these pairs, the enumerable information transfers do represent the plausible causal relationships.

**Proposition 4.2.** *There are strings $x, y$ and $x', y'$ for which the same $P$ is a minimal length conditional sufficient statistic such that*

$$
\begin{aligned}
EIT(x \leftarrow y) &=^+ 0 & EIT(x' \leftarrow y') &=^+ n - \log n \\
EIT(y \leftarrow x) &=^+ n - \log n & EIT(y' \leftarrow x') &=^+ 0.
\end{aligned}
$$

*Proof.* Let $k = \log n$ and let

$$
t_k = \min\{t : m(\epsilon, \epsilon) - m_t(\epsilon, \epsilon) \leqslant 2^{-k}\}.
$$

Remark that $t_k$ can be computed from $n$ and the first $k$ bits of $m(\epsilon, \epsilon)$, thus $K(t_k) \leqslant^+ k$. Let $a \in 2^n$ be the lexicographic string with $K_{t_k}(a) \geqslant n$. Let $x = x'$ such that $K(x|a^*) \geqslant^+ n$ and let

$$
\begin{aligned}
y &= \mathrm{XOR}(a, x_{2\ldots n}0) \\
y' &= \mathrm{XOR}(a, 0x^{n-1}).
\end{aligned}
$$

Than it follows that $K_{t_{k+O(1)}}(x, y) \leqslant^+ \frac{5}{4}n$. Moreover, by [4, Lemma 7.6]

$$
\begin{aligned}
K_{t_{k-O(1)}}(x, y) &\geqslant^+ K_{t_{k-O(1)}}(a, x) \\
&\geqslant^+ K_{t_k}(x|a) + K_{t_k}(a) \\
&\geqslant^+ 2n
\end{aligned}
$$

Since $K(x, y) \leqslant^+ \frac{5}{4}n$, this shows that the $m$-depth is $k$. By [3, Proposition ?.?] it follows that $P_{xy}$ as constructed there is a minimal sufficient statistic of $x, y$. The same reasoning shows that $P_{xy}$ is also a minimal sufficient statistic of $x', y'$.

It remains to show now the four inequalities of the Proposition. Remark that $K(x) =^+ n$, and also

$$
\begin{aligned}
K(y) &=^+ K(x, y) - K(x|y^*) \\
&\geqslant^+ K(x, a) - K(a) \\
&=^+ K(x|a^*) =^+ n,
\end{aligned}
$$

since also $K(y) \leqslant^+ l(y) = n$. Remark that

$$
\begin{aligned}
K(y|x \uparrow) &\leqslant^+ K(a) = k \\
K(x|y \uparrow) &\leqslant^+ l(x) = n,
\end{aligned}
$$

and since

$$
n + k \geqslant^+ K(x|y \uparrow) + K(y|x \uparrow) \geqslant^+ K(x, y) =^+ n + k,
$$

this shows that

$$
\begin{aligned}
K(y|x \uparrow) &=^+ k \\
K(x|y \uparrow) &=^+ n.
\end{aligned}
$$

A analogue argument shows the inequalities for $x', y'$. $\qquad\square$

The argument can be extended to show the analogue of Proposition 4.2 for the multivariate case with a complex incompressible causal structure, that contains no Halting information.

# 5 Conclusion

Ratio tests of universal semimeasures are defined and interpreted as an ideal solution to scientific composite hypothesis testing if an unlimited amount of computation can be assumed. Such universal elements exists if the hypotheses corresponds to a convex set or a product of convex sets of semimeasures, as for example in the hypothesis of independent sources for pairs of observations.

Using generalized structural equations, different hypotheses of influence and causality can be defined. These hypotheses define sets of enumerable semimeasures that have a universal element, and therefore they define different statistical tests. The detailed derivation of the statistical tests shows that there can be substantial differences in the corresponding confidences depending on the presumed directions of instantaneous information flows.

Associated causal semimeasures define a larger set of causal semimeasures who are not enumerable nor have a universal element. However, for the set of semimeasures associated with universal semimeasures, it is not clear whether a universal element exists, and consequently it is not clear whether they define some natural independence tests. However, these these tests can define ideal influence tests without assumptions on instantaneous information transfer. Different relations are summarized in figure 3.

Finally the ideal methods of information transfer are contrasted with practical methods from literature. Also the method is contrasted with the use of minimal sufficient statistics and it is shown that enumerable information transfer can describe plausible causal relations where minimal sufficient statistics can not do.

An upcoming paper provides some relatively tight coding theorem for these causal semimeasures which is exact within the logarithm of some notion of computational depth.

# Acknowledgments

The author was supported by a Ph.D grant of the Institute for the Promotion of Innovation through Science and Technology in Flanders (IWT-Vlaanderen). The author is grateful to Georges Otte, Luc Boullart and Patrick Santens for maintaining a marvelous interdisciplinary research group. I thank Bart Wyns and Erik Quaeghebeur for some linguistic improvements of the understanding of early versions of the manuscript.

# References

[1] B. Bauwens. Additivity of on-line decision complexity is violated by a linear term in the length of a binary string. *ArXiv e-prints*, August 2009.

[2] B. Bauwens. Influence tests II: $m$-depth and on-line coding results. In preparation, 2009.

[3] B. Bauwens. On the equivalence between minimal sufficient statistics, minimal typical models and initial segments of the Halting sequence. *ArXiv e-prints*, November 2009.

[4] B. Bauwens and S. Terwijn. Notes on sum-tests and independence tests. Accepted for publication in Theory of Computing Systems, 2009.

[5] B. Bauwens, B. Wyns, D. Devlaminck, G. Otte, L. Boullart, and P. Santens. Measuring instantaneous directed dependencies in interacting oscillators. In *Proceedings of the 28th Symposium on Information Theory in the Benelux*, 2007. http://www.autoctrl.ugent.be/bruno/papers/benelux.pdf.

[6] A. Chernov, S. Alexander, N. Vereshchagin, and V.Vovk. On-line probability, complexity and randomness. In *ALT '08: Proceedings of the 19th international conference on Algorithmic Learning Theory*, pages 138–153, Berlin, Heidelberg, 2008. Springer-Verlag.

[7] J. D. Collier. *Causation and Laws of Nature*, chapter Causation is the transfer of information, pages 215–246. Kluwer A. Publishers, 1999.

[8] M. Ding, Y. Chen, and S. L. Bressler. Granger causality: basic theory and application to neuroscience. In M. Winterhalder and J. Timmer, editors, *Handbook of Time Series Analysis*, pages 437–460. Wiley-cvh, Wienheim, 2006.

[9] M. Feder and N. Merhav. Universal composite hypothesis testing: a competitive minimax approach. *IEEE Transactions on information theory*, 48(6):1504–1517, 2002.

[10] U. Feldmann and J. Bhattacharya. Predictability improvement as an asymmetrical measure of interdependence in bivariate time series. *International Journal of Bifurcation and Chaos*, 14(2):505–514, 2004.

[11] P. Gacs. On the symmetry of algorithmic information. *Soviet Mathematical Dokledy*, 15:1477–1480, 1974.

[12] P. Gacs. Lecture notes on descriptional complexity and randomness. Unpublished, http://www.cs.bu.edu/faculty/gacs/papers/ait-notes.pdf, 2009.

[13] P. Gács, J. Tromp, and P.M.B. Vitányi. Algorithmic statistics. *IEEE Transactions on Information Theory*, 47(6):2443–2463, 2001.

[14] C.W.J. Granger. Investigating causal relations by econometric models and cross-spectral methods. *Econometrica*, 37:424–438, 1969.

[15] G. John. Inference and causality in economic time series models. In Z. Griliches and M. D. Intriligator, editors, *Handbook of Econometrics*, volume 2, chapter 19, pages 1101–1144. Elsevier, 1984.

[16] M. Kaminski, M. Ding, W.A. Truccolo, and S.L. Bressler. Evaluating causal relations in neural systems: granger causality, directed transfer function and statistical assessment of significance. *Biological Cybernetics*, 85(2):145–157, August 2001.

[17] R.J. Lawrence, S.A. Mulaik, and J.M. Brett. *Causal Analysis*. Sage Publications, 1983.

[18] J. Lemeire, K. Steenhaut, and A. Touhafi. When are graphical causal models not good models? In J. Williamson, F. Russo, and P. McKay, editors, *Causality in the Sciences*. Oxford University Press, 2010.

[19] L.A. Levin. Randomness conservation inequalities; information and independence in mathematical theories. *Inf. Control*, 61(1):15–37, 1984.

[20] M. Li and P.M.B. Vitányi. *An Introduction to Kolmogorov Complexity and Its Applications.* Springer-Verlag, New York, 2008.

[21] I. Martel. *Probabilistic Empiricism: In Defence of a Reichenbachian Theory of Causation and the Direction of Time.* PhD thesis, University of Colorado, 2000.

[22] J. Neyman and E. S. Pearson. On the Problem of the Most Efficient Tests of Statistical Hypotheses. *Royal Society of London Philosophical Transactions Series A*, 231:289–337, 1933.

[23] N.K.Verschagin. Agorithmic minimal sufficient statistics: a new definition. Presented on the 4th conference on randomness, computability and logic, Luminy., nov 2009.

[24] M. Palus and A. Stefanovska. Direction of coupling from phases of interacting oscillators: an information theoretic approach. *Physical Review E, Rapid Communications*, 67:055201(R), 2003.

[25] J. Pearl. *Causality: Models, Reasoning and Inference.* Cambridge University Press, 2000.

[26] M. G. Rosenblum and A. S. Pikovsky. Detecting direction of coupling in interacting oscillators. *Phys. Rev. E*, 64(4):045202, Sep 2001.

[27] K. Sameshima and L.A. Baccala. Using partial directed coherence to describe neuronal ensemble interactions. *Journal of Neuroscience Methods*, 94:93–103(11), Dec 1999.

[28] T. Schreiber. Measuring information transfer. *Physical Review Letters*, 85(2):461–464, Jul 2000.

[29] C. E. Shannon. A mathematical theory of communication, Jul and Oct 1948.

[30] D. A. Smirnov and B. P. Bezruchko. Estimation of interaction strength and direction from short and noisy time series. *Phys. Rev. E*, 68(4):046209, Oct 2003.

[31] P. Suppes. *A Probabilistic Theory of Causality.* North-Holland Publishing Company, 1970.

[32] M. Winterhalder, B. Schelter, W. Hesse, K. Schwab, L. Leistritz, R. Bauer, J. Timmer, and H. Witte. Comparison of linear signal processing techniques to infer directed interactions in multivariate neural systems. *Signal Processing*, 85:2137–2160, 2005.

## Prefix-free Kolmogorov complexity

For excellent introductions to Kolmogorov complexity we refer to [12, 20]. The definitions differ here, with respect that Kolmogorov complexities are conditioned on the parameter $n$, in most cases representing the length of the first argument.

An interpreter $\Phi$ is a partial computable function:

$$\Phi : \omega \times 2^{<\omega} \times \omega^{<\omega} \to \omega^{<\omega} : t, p, x \to \Phi_t(p|x).$$

and $\Phi(p|x) = \lim_{t\to\infty} \Phi_t(p|x)$. The use of $\omega^{<\omega}$ in this definition is to allow $\Phi$ to have multiple inputs and outputs in $\omega$ associated with $2^{<\omega}$. An interpreter is prefix-free if for any $x$, the set $D_x$ of all $p$ where $\Phi(p|x)$ is defined, is prefix-free. Among the prefix-free Turing machines, there are machines that can simulate any computation on any other prefix-free machine by prefixing $p$ with a finite binary string. These machines are called optimal universal Turing machines. Let $\Phi$ be some fixed optimal universal prefix-free interpreter.

For some $n \in \omega$, and $x, y \in \omega^{<\omega}$, the *Kolmogorov complexity* $K(x|y)$, is defined as:

$$
\begin{aligned}
K_t(x|y) &= \min\{l(p) : \Phi_t(p|y, n) \downarrow = x\} \\
K_t(x) &= K_t(x|\epsilon).
\end{aligned}
$$

$K(x|y)$ and $K(x)$ are obtained by taking the limit in $t$.

Some properties of length conditional prefix-free Kolmogorov complexity for $x \in 2^n$:

$$
K(x) \leqslant^+ n
$$
$$
K(x) =^+ K(K(x), x)
$$

Prefix-free Kolmogorov complexity is additive:

$$
K(x, y) =^+ K(x) + K(y|x^*), \tag{21}
$$

where $x^*$ is a program of length $K(x)$ that outputs $x$.

For $x, y \in \omega^{<\omega}$, $n \in \omega$,

$$
x \longrightarrow y
$$

means that there is a program $p_x$ with $l(p_x) \leqslant O(1)$, such that $\Phi(p_x|y, n) \downarrow = x$. Remark that $\Phi$ is also conditioned to $n$. Also remark that if $x \longrightarrow y$, than $K(x) \geqslant^+ K(y)$.

**Lemma 5.1.** *For any $w, p \in \omega$ with $\Phi(p) \downarrow = w$ and $l(p) \leqslant^+ K(w)$ we have*

$$
w^* \longleftrightarrow p \longleftrightarrow w, K(w). \tag{22}
$$

*Proof.* See [1, Lemma 4.2]. □

Let $m$ be a universal semimeasure The coding theorem states:

$$
-\log m(x|y) =^+ K(x|y).
$$